\documentclass[10pt,twoside]{article}
\usepackage{amssymb}
\usepackage{Latex-document}

\newtheorem{theo}{Theorem}[section]
\newtheorem{prop}{Proposition}[section]
\newtheorem{lemme}{Lemma}[section]

\renewcommand{\thesection}{\arabic{section}}
\newcommand{\Section}[1]{\section{\hskip -1em.\hskip 0.7em#1}}

\title{\bf Quasilinear Wave Equations and\vskip -2mm Microlocal Analysis\vskip 6mm}
\author{{\bf Hajer Bahouri}\vspace*{-0.5cm}\thanks{Universit\'e  \ de \ Tunis, \ D\'epartement \ de \
Math\'ematiques, \ 1060 \ Tunis, \ Tunisia. \ E-mail: Hajer.Bahouri@fst.rnu.tn}\vspace*{0.5cm} \quad Jean-Yves
Chemin\vspace*{-0.5cm}\thanks{Centre de Math\'ematiques \'Ecole polytechnique, 91 128 Palaiseau Cedex, France.
E-mail: chemin@math.polytechnique.fr}}
\date{\vspace{-8mm}}

\begin{document}

\def\NN{{\mathbb N}}
\def\RR{{\mathbb R}}
\def\EE{{\mathbb E}}
\def\HH{{\mathbb H}}
\def\CC{{\mathbb C}}
\def\ZZ{{\mathbb Z}}
\def\SS{{\mathbb S}}
\def\longformule#1#2{
\displaylines{ \qquad{#1} \hfill\cr \hfill {#2} \qquad\cr } }
\def\virgp{\raise 2pt\hbox{,}}
\def\cdotpv{\raise 2pt\hbox{;}}
\def\eqdefa{\buildrel\hbox{\footnotesize def}\over =}
\def\refer#1{~\ref{#1}}
\def\refeq#1{~(\ref{#1})}
\def\ccite#1{~\cite{#1}}
\def\sumetage#1#2{
\sum_{\scriptstyle {#1}\atop\scriptstyle {#2}} }
\def\Supp{\mathop{\rm Supp}\nolimits\ }
\def\suite#1#2#3{
({#1}_{#2})_{{#2}\in{#3}} }
\def\Id{{\rm Id}}
\def\ds{\displaystyle}

\let\al=\alpha
\let\b=\beta
\let\g=\gamma
\let\d=\delta
\let\ve=\varepsilon
\let\e=\epsilon
\let\z=\zeta
\let\lam=\lambda
\let\r=\rho
\let\s=\sigma
\let\t=\tau
\let\f=\phi
\let\vf=\varphi
\let\p=\psi
\let\om=\omega
\let\G= \Gamma
\let\D=\Delta
\let\Lam=\Lambda
\let\O=\Omega

\def\cA{{\mathcal A}}
\def\cB{{\mathcal B}}
\def\cC{{\mathcal C}}
\def\cD{{\mathcal D}}
\def\cE{{\mathcal E}}
\def\cF{{\mathcal F}}
\def\cG{{\mathcal G}}
\def\cH{{\mathcal H}}
\def\cI{{\mathcal I}}
\def\cJ{{\mathcal J}}
\def\cK{{\mathcal K}}
\def\cL{{\mathcal L}}
\def\cM{{\mathcal M}}
\def\cN{{\mathcal N}}
\def\cO{{\mathcal O}}
\def\cP{{\mathcal P}}
\def\cQ{{\mathcal Q}}
\def\cR{{\mathcal R}}
\def\cS{{\mathcal S}}
\def\cT{{\mathcal T}}
\def\cU{{\mathcal U}}
\def\cV{{\mathcal V}}
\def\cW{{\mathcal W}}
\def\cX{{\mathcal X}}
\def\cY{{\mathcal Y}}
\def\cZ{{\mathcal Z}}

\let\wt=\widetilde
\let\wh=\widehat
\let\tri\triangle
\def\ra{\rightarrow}
\def\pa{\partial}
\def\cdottri{\dot{\triangle}}
\def\cdotB{\dot{B}}
\def\cdotH{\dot{H}}
\def\cdotC{\dot{C}}
\def\mod#1{\langle{#1}\rangle}
\def\supetage#1#2{
\sup_{\scriptstyle {#1}\atop\scriptstyle {#2}} }

\def\cdv{champ de vecteurs}
\def\cdvs{champs de vecteurs}
\def\demo{d\'emons\-tra\-tion}

\maketitle

\thispagestyle{first} \setcounter{page}{141}

\begin{abstract}

\vskip 3mm

In this text, we shall give an outline of some recent results (see \ccite{bahourichemin2}  \ccite{bahourichemin3}
and \ccite{bahourichemin4}) of local wellposedness for two types of quasilinear wave equations for initial data
less regular than what is required by the energy method. To go below the regularity prescribed by the classical
theory of strictly hyperbolic equations, we have to use the particular properties of the wave equation. The result
concerning the first kind of equations must be understood as a Strichartz estimate for wave operators whose
coefficients are only Lipschitz while the result concerning the second type of equations is reduced to the proof
of a bilinear estimate for the product of two solutions for wave operators whose coefficients are not very
regular. The purpose of this talk is to emphasise the importance of ideas coming from microlocal analysis to prove
such results.

The method known to prove Strichartz estimates uses a representation eventually approximate of the solution. In
the case of the wave equation, the approximation used is the one coming from the Lax method, namely the one
connected to the geometrical optics. But it seems impossible, in the framework of the quasilinear wave equations,
to construct a suitable approximation of the solution on some interval $[0,T]$,  since the associate
Hamilton-Jacobi equation develop singularities (it is the caustic phenomenon) at a time connected with  the
frequency size. We have then to microlocalize, which means to localize in frequencies, and then to work on time
interval whose size depend on the frequency considered. It is the alliance of geometric optics and harmonic
analysis which allow to establish a quasilinear Strichartz estimate and to go below this minimal regularity in the
case of the first kind of equations.

To study the second kind of equations, we are confronted to an additional problem: Contrary to the constant case,
the support of the Fourier transform is not preserved by the flow of the variable coefficient wave equation. To
overcome this difficulty, we show that the relevant information in the variable case is  the concept of
microlocalized function due to J.M.Bony\ccite{bonyperso}. The proof that for solutions of variable coefficient
operators, microlocalization properties propagate nicely along the Hamiltonian flows related to the wave operator
is the key point in the proof of the result in the second case.

\vskip 4.5mm


\end{abstract}

\vskip 12mm

\Section{Introduction}

\vskip-5mm \hspace{5mm}

In this paper, our interest is to
prove local solvability for quasilinear wave equations of the type
\[ (E)\left\{
\begin{array}{rcl}
\partial^2_t u -\D u -g(u)\cdot\nabla^2 u & = & Q(\partial u,\partial u)\\
(u,\partial_t u) _{|t=0} & = & (u_0,u_1)
\end{array}
\right.
\]
where~$g$ is a smooth function  vanishing at~$0$ with value in~$K$
such that~$\Id+K$ is a convex compact subset of the set of
positive symmetric matrices and ~$Q$ is a  quadratic form
on~$\RR^{d+1}$. Our interest proceeds also for cubic quasilinear
wave equations of the type
\[(EC)
\left\{
\begin{array}{rcl}
\ds
\partial_t^2u -\D u -\sum_{1\leq j,k \leq d} g^{j,k}\partial_j\partial_k u
& = & \ds \sum_{1\leq j,k \leq d} \wt Q_{j,k} (\partial g^{j,k},
\partial u) \\ \D g^{j,k} & = & Q_{j,k}(\partial u,\partial u) \\
(u,\partial_t u)_{|t=0}& = & (u_0,u_1)\\
\end{array}
\right.
\]
where~$Q_{j,k}$ and ~$\wt Q_{j,k}$  are quadratic forms
on~$\RR^{d+1}$ and where all the quadratic forms are supposed to
be smooth functions of~$u$. \smallbreak The basic tool to prove
local solvability for such equations is the following energy
estimate, also valid for the symmetric systems
\begin{equation}
\|\partial u (t, . )\|_{ H^{s-1}} \leq \|\partial u (0, . )\|_{
H^{s-1}} e^{C \int_0^T \|\partial g (\tau, .
)\|_{L^{\infty}}d\tau . }
\end{equation}
So thanks to classical arguments, local solvability derives easily
from the control of the quantity \[ \int_0^T \|\partial g (\tau, .
)\|_{L^{\infty}}d\tau .
\]
In the framework of the equation ~$(E)$,  the control of this key
quantity  requires initial data~$(u_0,u_1)$ in~$H^{s}\times
H^{s-1}$ for ~$ s > \frac d 2+1 $ while in the framework of~$(EC)$
(with small data, which makes sense in this case) it only requires
initial data~$(u_0,u_1)$ in~$H^{\frac d 2+\frac 1 2}\times
H^{\frac d 2- \frac 1 2}$.

The goal of this paper is to go below this regularity for the
initial data. Let us first have a look at the scaling properties
of equations~$(E)$ and ~$(EC)$. If~$u$ is a solution of~$(E)$
or~$(EC)$, then~$u_\lam(t,x)\eqdefa u(\lam t,\lam x)$ is also a
solution of the same equation. The space which is invariant under
this scaling for the couple ~$(u_0,u_1)$ is ~$\dot H^{\frac d
2}\times \dot H^{\frac {d}{ 2}-1}$. So the results given by the
classical energy estimate appear to require more regularity than
the scaling in the two cases.

In fact, the energy methods despise the particular properties of
the wave equation. It is on the impulse of the pioneer work of S.
Klainerman (see\ccite{klainerman1}) that a vast series of works
have been attached to improve the span life time of regular
solutions of quasilinear wave equations using the Lorentz
invariance. Let us notice the results of S. Alinhac
(see\ccite{alinhac11} and\ccite{alinhac12}), of
 L. H\"ormander (see\ccite{hormander2}),  of F. John (see\ccite{john2}),
  of F. John and S. Klainerman
(see\ccite{johnklainerman}), of S. Klainerman
(see\ccite{klainerman2}) and of J-M.Delort (\ccite{delort1})
concerning the Klein-Gordon equation.

In this talk, we shall limit our self to the question of minimal regularity. Concerning this subject, the only
case studied is the semilinear case,  which means the case of the equation~$(E)$ with~$g\equiv 0$. As it has been
shown by S. Klainerman and M.Machedon (see\ccite{klainermanmachedon3} and \cite{klainermanmachedon5}) we can, when
the quadratic form ~$Q$ verifies a  structure condition known by ``null condition", nearly reach the space
invariant by scaling.  For any quadratic form~$Q$, we have the following theorem, proved by G. Ponce and T.
Sideris in\ccite{poncesideris}
\begin{theo}
\label{ondesemilinstritheo} Let us define~$\underline s_d$ by  $$
\underline s_{\,d}\eqdefa  \frac d 2+\frac 1 2 \quad\hbox{if}\quad
d\geq 3 \quad\hbox{and}\quad \underline s_{\,2}= \frac 7 4\cdotp
$$ Let~$(u_0,u_1)$ be a Cauchy data in ~$H^s\times H^{s-1}$
with~$s>\underline s_{\,d}$ then there exists a time~$T$ such that
there exists a unique  solution~$u$ of the equation~$(E)$ such
that $$ u\in L^\infty([0,T];H^s)\cap
Lip([0,T];H^{s-1})\quad\hbox{and}\quad \partial u\in
L^2([0,T];L^\infty). $$
\end{theo}

The proof of this result lies on specific properties of the wave
equation, namely the following Strichartz estimate
\begin{eqnarray}
  \|\partial u \|_{L^2_T(L^{\infty})} &\leq &C \left( \|\partial u(0,.) \|_{
  H^{s-1}} + \|\Box\, u \|_{L^1_T(H^{s-1})}\right), \nonumber \\
  & & \hspace*{2cm} \mbox{for $d \geq 3$ and $s > \frac d 2 + \frac 1 2$}.
\end{eqnarray}
Indeed, if we couple it with the standard energy estimate
 \[\|\partial u \|_{L^{\infty}_T(H^{s-1})} \leq C \left( \|\partial u(0,.) \|_{
 H^{s-1}} + \|\Box\, u \|_{L^1_T(H^{s-1})}\right)
 \]
 we obtain, owing to the tame estimates and the Cauchy-Schwarz inequality
\begin{eqnarray*}
 & &\|\partial u \|_{L^2_T(L^{\infty})} + \|\partial u
 \|_{L^{\infty}_T(H^{s-1})}\\
 &\leq &C \left( \|\partial u(0,.) \|_{
 H^{s-1}} + T^{\frac 1 2} \|\partial u(t,.) \|_{L^2_T(L^{\infty})}\|\partial u
 \|_{L_T^{\infty}H^{s-1}}\right),
\end{eqnarray*}
 which ensures by the theory of evolution equations the local
 solvability for ~$T \leq \frac{C}{\|\partial u(0,.) \|_{
 H^{s-1}}^2}.$

In other respects, in\ccite{lindblad1}, H. Linblad shows that
for~$d=3$ the above result is optimum,  which means that the
problem~$(E)$ with~$g\equiv 0$ is  not wellposed in ~$H^2$. Let
us also notice that the same kind of result is also true on the
Heisenberg group (see \ccite{bahourigallagher}).

The authors (see\ccite{bahourichemin2} and\ccite{bahourichemin3})
adjust a method followed by D. Tataru (see\ccite{tataruondeql2})
based on microlocal analysis to improve the minimal regularity
for the equation~$(E)$ in the quasilinear case. Let us recall
this result
\begin{theo}
\label{ondequasilinstritheo} {\it If~$d\geq 3$, let~$(u_0,u_1)$ be
in~$H^s\times H^{s-1}$ for~$s>s_d$ with~$s_d = \displaystyle \frac
d 2+\frac 1 2+\frac 1 6\,\cdotp$ Then, a positive time~$T$ exists
such that a unique solution~$u$ of  the equation~$(E)$ exists such
that
\[
\partial u\in C([0,T];H^{s-1})\cap  L^2([0,T];L^\infty).
\]
 }
\end{theo}
\noindent{\bf Remarks }
\begin{itemize}
\item This theorem has been proved with ~$1/4$ instead
than~$1/6$ in\ccite{bahourichemin2} and then improved a little bit
in\ccite{bahourichemin3} and proved with~$1/6$ by D. Tataru
in\ccite{tataruondeql2}. Strichartz estimates for quasilinear
equations are the key point of the proofs.
\item Let us notice that the improvement with~$1/6$ of  D. Tataru
in\ccite{tataruondeql2} is  due to a different manner of counting
 the intervals where microlocal estimates are true.
\item Recently, S. Klainerman and I. Rodnianski in\ccite{klainermanrodnianski} have
obtained a better index in dimension 3. Their proof is based on
very different methods.
\item Let us notice that we have also improved the minimal regularity in dimension
2, but the gain is only of~$\frac 1 8 $ derivative, this is
explained by the mean dispersif effect in this dimension already
known for the constant case.
\end{itemize}
 The analogous theorem in the case
of equation~$(EC)$ is the following
\begin{theo}
\label{ondecubiquesimpledequal4} {\it If~$d\geq 4$,
let~$(u_0,u_1)$ be in~$H^{s}\times H^{s-1}$ with~$ s > \frac d 2
+\frac 1 6$ such that~$\|\g\|_{\dot H^{\frac d 2 -1}}$ is small
enough. Then, a positive times~$T$ exists such that a unique
solution~$u$ of~$(EC)$ exists such that\[\partial u \in
C([0,T];H^{s-1})\cap L^2_T(\dot B^{\frac d 4-\frac 1 2}_{4,2})
,\quad for \quad  d\geq 5,
\]
and
\[
\partial u \in C([0,T];H^{s-1})\cap L^2_T(\dot B^{\frac 1 6}_{6,2})
\quad\hbox{and}\quad \partial g \in L^1_T(L^\infty) \quad for
\quad d=4.
\]
 }
\end{theo}
\noindent{\bf Remarks }
\begin{itemize}
\item
The case when~$d\geq 5$ can be treated  only with Strichartz
estimates simply because if~$\partial u$ belongs to~$L_T^2(\dot
B^{\frac d 4-\frac 1 2}_{4,2})$ then~$\partial g$ is
in~$L^1_T(L^\infty)$.
\item
The case when~$d=4$ requires bilinear estimates. This fact appears
in the statement of Theorem\refer{ondecubiquesimpledequal4}
through the following phenomenon: The fact that~$\partial u$ is
in~$L^2_T(\dot B^{\frac 1 6}_{6,2})$ does not imply that the time
derivative of~$g$ belongs to~$L^1_T(L^\infty)$. Of course this
condition is crucial in particular to get the basic energy
estimate. But we have been unable to exhibit a Banach space~$\cB$
which contains the solution~$u$ and such that if a function~$a$ is
in~$\cB$, then~$\partial \D^{-1}(a^2)$ belongs
to~$L^1_T(L^\infty)$.
\item For technical obstructions, this theorem is limited to the
dimensions~$d \geq 4$.
\end{itemize}

\Section{Quasilinear Strichartz estimates}

\vskip-5mm \hspace{5mm}

Following the process of  G. Ponce and T. Sideris
in\ccite{poncesideris}, we reduce the proof of the theorem
\refer{ondequasilinstritheo} to the following a priori estimate
\begin{theo}
\label{quasistrichartz}
 {\it If~$d\geq 3$, a constant C  exists such that, for any regular solution u of
  the equation ~$(E)$, if
  \[ T^{\frac 1 2 +(s-s_d)}\left(\|\g\|_{ H^{\frac d 2 -\frac 1 2 +
  (s-s_d)}}+ T^{\frac 1 6}
 \|\g\|_{ H^{s-1}}\right) \leq C, \quad with \quad s > s_d
  \]
  then we have \[\|\partial u \|_{L^2_T(L^{\infty})}
 \leq C \left( \|\g\|_{H^{s-1}} + \|Q(\partial u,\partial u) \|_{L^1_T(H^{s-1})}\right).\]
   }
   \end{theo}

The estimate in hand must be understood as a Strichartz estimate
for wave equations with variable coefficients and not very
regular. The Strichartz estimates have a long history begining
with  Segal's work \ccite{segal} for the wave equation with
constant coefficients. After the fundamental work of Strichartz
\ccite{strichartz}, it was developed by diverse authors, we refer
to the synthesis article of Ginibre and Velo \ccite{ ginibrevelo2}
to which it is advisable to add the recents works of Keel-Tao
\ccite{KeelTao} consecrated to some limited cases and of Bahouri,
G\'erard and Xu \ccite{bgx} for the wave equation on the
Heisenberg group. For Strichartz estimates with ~$C^{\infty} $
coefficients, we refer to the result of L. Kapitanski
(see\ccite{kapitanski}). The article of H. Smith
(see\ccite{smith1}) constitutes an important step in the study of
Strichartz estimates for  operators with coefficients not very
regular since it proves  Strichartz estimates with coefficients
only~$C^{1,1}$.

We shall now explain how to establish this quasilinear Strichartz
estimate, showing where are the difficulties and what are  the
essential ideas which allow to overcome them. The method known to
prove these estimates uses a representation, eventually
approximate, but always explicit of the solution. In the case of
the wave equation, the approximation used is the one coming from
the Lax method, namely the one connected to the geometrical
optics. To make such a method work in the framework of
quasilinear wave equations requires a ``regularization'' of the
coefficients also in time. This leads to the following iterative
scheme introduced in\ccite{bahourichemin3}.
 Let us define the sequence~$ (u^{(n)})_{n\in \NN}$ by the first term~$u^{(0)}$
satisfying
\[
\left\{
\begin{array}{rcl}
\partial_t^2 u^{(0)}-\D u^{(0)} & = & 0\\
( u^{(0)}, \partial_t u^{(0)})_{|t=0} & = & (S_0u_0,S_0u_1),
\end{array}
\right.
\]
and by the following induction
\[
(E_n)\,\left\{\begin{array}{rcl}
\partial_t^2 u^{(n+1)} -\D u^{(n+1)}
-g_{n,T}\cdot\nabla^2 u^{(n+1)}  & = & 0\\
 (u^{(n+1)},\partial_tu^{(n+1)})_{|t=0} & = & (S_{n+1}u_0,S_{n+1}u_1)
\end{array}
\right.
\]
where~$g_{n,T}\eqdefa\theta(T^{-1})g_n$ with~$g_n \eqdefa g(u^n)$
and~$\theta$ a function of~$\cD(]-1,1[)$ whose value is~$1$
near~$0$ and where ~$S_n $ is a frequencies truncated  operator
which only  conserves the frequencies lower than ~$ C 2^{n-1}$.
Let us introduce some notations which will be used all along this
section. If ~$s=s_d + \al$ where~$\al$ is a small positive number,
let us define
\[
 N_T^\al(\g) \eqdefa \|\g\|_{H^{\frac d 2 -\frac 1 2 +\al}}+ T^{\frac 1 6 } \|\g\|_{H^{s-1}}.
\]

The assertion we have to prove by induction, for ~$T^{\frac 1 2
+\al }N_T^\al(\g)$ small enough, is the following: If~$ d\geq 3$,
\[
(\cP_n)\left\{
\begin{array}{ccl}
\|\partial u^{(n)}\|_{L^2_T (L^{\infty })} & \leq & C_{\al}T^{\al}
N_T^\al(\g)\\ \| \partial u^{(n)}\|_{L^{\infty}_T (H^{s-1})} &
\leq & C \|\g\|_{H^{s-1}.}
\end{array}
\right.
\]

For this, we shall transform the equation~$(E_n)$ into a
paradifferential equation, more precisely an equation of the type
\[
(E_q)\qquad \partial^2_t u^{(n+1)}_q - \Delta u^{(n+1)}_q -
\left(S_q g_{n,T}\right)\nabla^2 u_q^{(n+1)} = \wt R_q(n)
\]
where the term $\widetilde R_q(n)$ is a remainder term estimated
as agreed and  where~$u_q$ denotes the part of~$u$ which is
relative to the frequencies of size  ~$2^q$.

This transformation of the equation, which is the classical
  paralinearization  defined by
J.-M. Bony in\ccite{bony2} is here not sufficient, since it is
well known that the paradifferential operators defined
in\ccite{bony2} belong to a bad class of
 pseudodifferential operators ( class $S_{1,1}$ of
H\"ormander), class in particular devoid of any asymptotic calculus, which forbids of course to envisage any
approximate method of type ``Lax method".

The idea is as in\ccite {lebeau3} to truncate more in the
frequencies of the metric~$g$ and to transform the equation
$(E_q)$ on the following equation $(EPM_q)$
\[
(EPM_q)\qquad \partial^2_t u^{(n+1)}_q - \Delta u^{(n+1)}_q -
\left(S_{ \delta q}g_{n,T}\right)\nabla^2 u^{(n+1)}_q = R_q(n);
\]
where ~$ 0 < \delta <1 $  and ~$S_{ \delta q}$ is a frequencies
truncated operator  which conserves only the frequencies smaller
than~$C T^{-(1-\delta )}2^{ \delta q -1}$. We can  interpret it as
a localization in the pseudodifferential calculus sense ~$(1,
\delta )$ of H\"ormander.

This localization allows us to  construct an approximation of the
 solution but engenders a loss in the remainder~$ R_q(n)$. This
 approximation is on the form $$ \int e^{i\Phi_q(t,x,\xi)} \sigma_q(t,x,\xi)
\widehat\gamma(\xi)d\xi $$ where $\Phi_q$ is a solution of the
Hamilton-Jacobi equation and $\sigma_q$ is a symbol calculated by
resolving a sequence of transport equations; it is about a
classical method. But, on account of the the caustic phenomenon,
this approximation is microlocal, which means  valid only a time
interval whose length depends on the size of the frequencies we
work with.

Nevertheless, following the classical method, we prove microlocal
Strichartz estimates
\begin{equation}
\label{microstrichartz} \|\partial u^{(n+1)}_q
\|_{L^2_{I_q}(L^\infty)} \leq C_{\beta}(2^q T)^{\beta} 2^{q (\frac
{d-1}{2})}\left( \|\g_q \|_{L^2} + \|R_q(n)
\|_{L^1_{I_q}(L^2)}\right)
\end{equation}
for any positive~$\beta$, where ~$\g_q \eqdefa (\nabla
(u_0)_q,(u_1)_q)$ and ~$I_q $ satisfies
\begin{equation}
 \label{conditionHJ}
 \|\nabla ^2 G^{(n)}_\delta \|_{L^1_{I_q}(L^\infty)} \leq \e
  \end{equation}
  where ~$G^{(n)}_\delta \eqdefa S_{ \delta q}g(u^{(n)})$ and
\begin{equation}
 \label{conditiondecoupst0}
|I_q | \leq T (2^q T)^{1-2\delta -\e}.
 \end{equation}
The condition\refeq{conditionHJ} is imposed by the
 Hamilton-Jacobi equation while the
 condition\refeq{conditiondecoupst0} is required by the asymptotic
 calculus to turn out the ``Lax method".

Finally to prove the complete estimate, the method we used consists in a decomposition of the interval~$[0,T]$ on
subintervals~$I_q$ on which the above microlocalized estimates are true. The key point is a careful counting  of
the number of such intervals, for this we shall use here D. Tataru's version of the method we introduced
in\ccite{bahourichemin2}.

The idea consists to seize at the opportunity of this decomposition to compensate the loss on the remainder. To do
so, we impose on the interval ~$I_q$ the supplementary condition
\begin{equation}
\|R_q (n)\|_{L^1_{I_q}(L^2)} \leq \lam \|R_q (n)\|_{L^1_T(L^2)}
\end{equation}
where the parameter~$\lam$ is to be determined in the
interval~$[0,1]$. This constraint joint to the
conditions\refeq{conditionHJ} and\refeq{conditiondecoupst0} leads
by optimization to the best choice
\[
\lam = (2^qT)^{-\frac \d 2}, \quad\quad \d = \frac 2 3\, ,
\]
and allows to conclude that the number~$N$ of such intervals is
less than~$C(2^qT)^{\frac 1 3 +\e}$. If we denote
by~$(I_{q,\ell})_{1\leq \ell \leq N}$ the partition of the
interval~$[0,T]$ on such intervals, we can write thanks to
\refeq{microstrichartz},
\begin{eqnarray*}
 \|\partial u^{(n+1)}_q  \|^2_{L^{2}_T (L^{ \infty})} &
 \leq & C_{\beta} \sum_{\ell =1}^{N}(2^q T)^{2\beta} 2 ^{2 q( \frac
 {d-1}{2})}  \left( \| \g_q \|_{L^{2}} + \| R_q(n) \|_{L^{1}_{I_{q,\ell}}
 (L^{2})}\right)^2
 \\&
 \leq & C_{\beta} \sum_{\ell =1}^{N} (2^q T)^{2\beta} 2 ^{2 q( \frac
 {d-1}{2})} \left( \| \g_q \|_{L^{2}} + (2^q T )^{- \frac 1 3 }
  \| R_q(n) \|_{L^{1}_T
 (L^{2})}\right)^2 .
\end{eqnarray*}
As ~$ N $ is less than ~$C(2^q T)^{\frac 1 3
 +\e}$, we obtain
 \begin{eqnarray}
 \|\partial u^{(n+1)}_q  \|_ {L^{2}_T (L^{ \infty})}
 &\leq & C_{\beta} (2^q T)^{\beta} 2 ^{ q( \frac
 {d-1}{2})}(2^q T )^{\frac 1 6
 + \frac {\e}{2}} \nonumber \\
 & &\cdot \left( \| \g_q \|_{L^{2}} + (2^q T )^{- \frac 1 3 }
  \| R_q(n) \|_{L^{1}_T
 (L^{2})}\right).
 \end{eqnarray}

 Now as the loss in the remainder ~$R_q(n)$ is of
 order ~$2^{q (1- \delta )}$ and more precisely, for ~$N^{\al}_T(\g)$ small enough,
  we have
 \begin{eqnarray}
  \| R_q(n) \|_{L^{1}_T(L^{2})} &\leq & c_q C 2^{-q(\frac
  {d-1}{2})} (2^q T)^{-(s - \frac {d}{2}- \frac {3}{2}+\delta )}
T^{s-\frac{d}{2}-\frac{1}{2}}   \| \g \|_{H^{s-1}} \nonumber \\
& &\cdot \left( 1 +T^{\frac 1 2 } \| \partial
u^{n+1}\|_{L^{2}_T(L^{\infty})}\right)
 \end{eqnarray}
 where~$(c_q) \in \ell^2$.
 We deduce, owed to the choice of ~$\delta$ that
 \[  \|\partial u^{(n+1)}_q  \|_{L^{2}_T (L^{ \infty})} \leq c_q C  ( 2
 ^{q}T)^{-(\al - \frac {\e}{2}-\beta)}
  T^{s-\frac{d}{2}-\frac{1}{2}} \| \g \|_{H^{s-1}} \left(1+ T^{\frac 1 2 }
 \| \partial u^{n+1} \|_{L^{2}_T(L^{\infty})}\right)\]
 which implies the result by summation.

\Section{Quasilinear bilinear estimates}

\vskip-5mm \hspace{5mm}

 The method used here is not without any interaction with the one used
to prove the theorem \refer{ondequasilinstritheo}. As in the case
of equation~$(E)$, the basic fact is the control of
\[
\int_0^T\|\partial g(\tau,\cdot) \|_{L^\infty} d\tau,
\]
and the proof of the theorem \refer{ondecubiquesimpledequal4}
follows from the following a priori estimate
\begin{theo}
\label{quasibilin} {\it If~$d\geq 4$, a constant C exists such
that, for any regular solution u of
  the equation ~$(EC)$, if ~$\|\g\|_{\dot H^{\frac d 2 -1}}$ is small
enough and
  \[ T^{\frac 1 6 +(s-\frac d 2
-\frac 1 6)}
 \|\g\|_{H^{s-1}} \leq C, \quad with \quad s >\frac d 2
+\frac 1 6
  \]
  then we have \[\|\partial \D^{-1} Q(\partial u , \partial u) \|_{L^1_T(L^{\infty})}
 \leq C  \|\g\|^2_{ H^{s-1} .} \]
   }
   \end{theo}

This is the quasilinear version of the following  bilinear
estimate
   owed to D. Tataru and S. Klainerman
(see\ccite{klainermantataru})
\begin{prop}
\label{klainermantatarubilin} {\it Let~$u$ be a solution
of~$\partial_t^2u -\D u =  0$ and~$(\partial u)_{|t=0} = \g$.
Then, if~$d\geq 4$,
\[
\|\partial \D^{-1} Q(\partial u , \partial u)\|_{L^1_T(L^\infty)}
\leq C_{\e,T}\|\g\|^{2}_{\frac d 2-1+\e} .
\]
}
\end{prop}
\noindent{\bf Remark } We find a gain of one  derivative from the
regularity of the initial data compared with the product laws and
a gain of half a derivative about the regularity of the initial
data compared with purely Strichartz methods.

To explain the basic ideas of bilinear estimates, let us first consider the case of constant coefficients.
As~$\partial_t \Delta^{-1} \bigl(\partial_j u(t,\cdot) \partial_k u(t,\cdot)\bigr) = \Delta^{-1}
\bigl(\partial_t\partial_j u\partial_k u(t,\cdot) \bigr) +$ \linebreak $\Delta^{-1} \bigl(\partial_j u
\partial_t\partial_k u(t,\cdot) \bigr)$, we have to control
expression of the type
\[
\int_0^T \|\Delta^{-1} \bigl(\partial_t\partial_j u\partial_k
u(t,\cdot) \bigr)\|_{L^\infty} dt.
\]
For this we introduce Bony's decomposition which consists in
writing
\[
ab = \sum_q S_{q-1}  a \D_q b +\sum_q S_{q-1}  b \D_q a +
\sumetage {-1\leq j\leq 1}{q} \D_qa \D_{q-j}b.
\]
 When~$d\geq 4$, we have~$ \|\partial^k u_q
\|_{L^2_T(L^\infty)} \leq C 2^{q\left(\frac d 2 -\frac 1
2+k-1\right)}\|\g_q\|_{L^2}$, then it is easy to  prove that
\[
\Bigl\|\D^{-1}\Bigl(\sum_q S_{q-1}  \partial^2 u  \partial
u_q\Bigl)\Bigr\|_{L^1_T(L^\infty)} \leq C\|\g\|_{\frac d 2-1}^2.
\]
The symmetric term can be treated exactly along the same lines.
The  remainder term
\begin{equation}
\label{bilinearremainder} \D^{-1}\Bigl(\sumetage {-1\leq j\leq
1}{q}
\partial^2 u_q
\partial u_{q-j} \Bigl)
\end{equation} is much more difficult to treat in particular in
dimension~$4$. The idea introduced by D. Tataru and S. Klainerman
(see\ccite{klainermantataru}) consists to treat this term using
precised Strichartz estimates and interaction lemma. \smallbreak
The precised Strichartz estimates are described by the following
proposition.
\begin{prop}
\label{corstrichatz1+} {\it Let~$\cC$ be a ring of~$\RR^d$.
If~$d\geq 3$, a  constant~$C$ exists such that for any~$T$ and
any~$h\leq 1$, if ~$\Supp \wh u_j$  are included in a ball of
radius~$h$ and in the ring~$\cC$, we have $$
\|u\|_{L^2_T(L^\infty)}\leq  C \bigl(h^{d-2}\log
(e+T)\bigr)^{\frac 1 2} \bigl(\|u_0\|_{L^2}+\|u_1\|_{L^2}\bigr),
$$where~$u$ denotes the solution of ~$\partial_t^2 u -\D u =0$ and
~$\partial_t^ju_{|t=0}=u_j$. }
\end{prop}

As usual it is deduced with the~$TT^\star$ argument from the
following dispersive inequality.
\begin{lemme}
\label{lemmestrichatz1+} {\it A constant~$C$ exists such that
if~$u_0$ and~$u_1$ are functions in~$L^1(\RR^d)$ such that
\[
\Supp (\wh u_j)\subset \cC\quad\hbox{and}\quad \max
\bigl\{\d(\Supp(\wh u_0)),\d(\Supp( \wh u_1))\bigr\}\leq h ,
\]
then, for any~$\wt d$ between~$0$ and~$d-1$, we have
\[
\|u(t,\cdot)\|_{L^\infty}\leq \frac {Ch^{d-\wt d}} {t^{\frac {\wt
d}2} } \bigl(\|u_0\|_{L^1}+\|u_1\|_{L^1}\bigr),
\]
where~$u$ denotes the solution of ~$\partial_t^2 u -\D u =0$ and
~$\partial_t^ju_{|t=0}=u_j$. }
\end{lemme}

This inequality is proved in\ccite{klainermantataru} in the case
~$\wt d=d-1$.  The general case is obtained by interpolation with
the classical Sobolev embedding.

Let us now show how to take account the interactions of the solutions to control the accumulation of frequencies
at the origin in the study of the remainder term.
\begin{lemme}``Interaction Lemma "
 {\it There exists a constant C such that if ~$v_1$ and~$v_2$ are two
solutions of~$\partial_t^2v_j -\D v_j =  0$ satisfying~$(\partial
v_j)_{|t=0} = \g_j$ with ~$\Supp (\wh {\g_j})\subset \cC$ we have
for, ~$ 0 < h < 1$
\begin{equation}
\label{interactioninproduct} \|\chi(h^{-1}D) (\partial^2v_1
\partial v_2)\|_{L^1_T(L^\infty)} \leq C h^{d-2}\log(e+T)
\|\g_1\|_{L^2} \|\g_2\|_{L^2},
\end{equation} where ~$\chi $ is a radial function in ~$\cD $
which is equal to 1 near the origin.}
\end{lemme}

Let us define~$(\f_\nu)_{1\leq \nu\leq N_h}$ a partition of unity
of the ring~$\cC$ such that~$\Supp \f_\nu \subset B(\xi_\nu,h)$.
Then, using the fact that the support of the Fourier transform of
the product of two functions is included in the sum of the
support of their Fourier transform, a family of
functions~$(\wt\f_\nu)_{1\leq \nu\leq N_h}$  exists such
that~$\Supp \wt \f_\nu \subset B(-\xi_\nu,2h)$ and
\begin{equation}
\label{interactprodcoefconstant}\chi(h^{-1}D) (\partial^2v_1
\partial v_2) = \sum_{\nu=1}^{N_h}\chi(h^{-1}D) \bigl(\partial^2
\wt \f_\nu(D)v_1 \partial \f_\nu(D)v_2\bigr).
\end{equation}
Applying Proposition\refer{corstrichatz1+} gives
\[
\|\chi(h^{-1}D) (\partial^2v_1 \partial v_2)\|_{L^1_T(L^\infty)}
\leq C h^{d-2}\log(e+T) \sum_{\nu=1}^{N_h}
\|\wt\f_\nu(D)\g_1\|_{L^2}\|\f_\nu(D)\g_2\|_{L^2}.
\]
The Cauchy Schwarz inequality implies that
\begin{eqnarray*}
  & &\|\chi(h^{-1}D) (\partial^2v_1 \partial v_2)\|_{L^1_T(L^\infty)}\\
  &\leq & C h^{d-2}\log (e+T) \Biggl(\sum_{\nu=1}^{N_h}
  \|\wt\f_\nu(D)\g_1\|^2_{L^2}\Biggr)^{\frac1
  2}\Biggl(\sum_{\nu=1}^{N_h}
  \|\f_\nu(D)\g_2\|^2_{L^2}\Biggr)^{\frac 1 2}.
\end{eqnarray*}
The almost orthogonality of~$(\wt \f_\nu(D)\g_1)_{1\leq \nu\leq N_h}$ and~$(\f_\nu(D)\g_2)_{1\leq \nu\leq N_h}$
implies \refeq{interactprodcoefconstant} and leads then to the estimate of the remainder
term\refeq{interactioninproduct} by rescaling.

To establish the theorem\refer{quasibilin}, we shall follow the
steps of the proof of the theorem\refer{quasistrichartz} which
consists owed to the gluing method to reduce the problem to the
proof of ``microlocal" bilinear estimates. The generalization of
the precised Strichartz estimates to the framework of the
equation~$(EC)$ doesn't cost more than the generalization of the
Strichartz estimates to the framework of the equation~$(E)$, the
supplementary difficulty to study the equation ~$(EC)$ lies in
the generalization of the interaction lemma. The preservation of
the support of the Fourier transform by the flow of the wave
equation is the crucial point in the proof of this lemma. The
defect of this property in the case of the variable coefficients
constitutes the additional major problem in the proof of the
theorem \refer{quasibilin}.

To palliate this difficulty, we have used a finer localization in
phase space. This localization is given by the concept of
microlocalized function near a point~$X=(x,\xi)$ of the cotangent
space~$T^\star \RR^d$ (the cotangent space of~$\RR^d$). More
precisely, if we consider the positive quadratic form g on
~$T^\star \RR^d$ defined by
\[
g(dy^2,d\eta^2) \,\eqdefa\, \frac {dy^2} {K^2} +\frac
{d\eta^2}{h^2} \quad\hbox{with}\quad \lam \eqdefa K h \geq 1
\]
a function~$u$ in~$L^2(\RR^d)$ is said to be microlocalized
in~$X_0 =(x_0,\xi_0)$ a point of~$T^\star \RR^d$ if
\[
\cM^{C_0,r}_{X_0,N}(u) \eqdefa \sup_{g(X-X_0)^{\frac 1 2} \geq
C_0r}\lam^{2N}g(X-X_0)^N\supetage{\vf\in
\cD(B_g(X,r))}{\|\vf\|_{k_N,g}\leq 1} \|\vf^Du\|_{L^2}
\]
are finite, where $B_g(X,r)$ denotes the  g-ball of center X and
radius r, the operator~$\vf^D$ is defined by
\[
(\vf^D u) (x) = (2\pi)^{-d}
\int_{T^\star\RR^d}e^{i(x-y|\xi)}\vf(y,\xi) u(y)dyd\xi,
\]
and
\[
\|\vf\|_{j,g} \eqdefa \supetage{k\leq j}{X\in
T^\star\RR^d}\supetage{(T_\ell)_{1\leq \ell\leq k}}{g(T_\ell)\leq
1} |D^k\vf(X)(T_1,\cdot,T_k)|.
\]

This notion  due to J.-M.Bony~(\cite{bonyperso}) means
 that the function ~u is concentrated in space near the point~$x_0$
 and in frequency near the point~$\xi_0$ and behaves well
  against the product, namely, we show that if \[
g(\check Y_1-Y_2)^{\frac 1 2}\geq C_0r,
\]  where~$\check Y
\eqdefa (y,-\eta)$ if~$Y=(y,\eta)$ then  for any~$N$, we have
\[\begin{array}{rcl}
& & \bigl\|\chi(h^{-1}D)(\vf^D_1 u_1\vf_2^D u_2)\bigr\|_{L^1}\\
&\leq & C_N\|\vf_1\|_{k_N,g}\|\vf_2\|_{k_N,g}
\bigl( 1+\lam^2 g(\check
Y_1-Y_2)\bigr)^{-N}\|u_1\|_{L^2}\|u_2\|_{L^2}\end{array}\] where
~$\vf_i \in \cD(B_g(Y_i,r))$.

This study of the interaction between two typical examples of
microlocalized functions allows as
in\refeq{interactprodcoefconstant} to concentrate the bilinear
estimate on real interaction.

Anyway, the choice of the localization metric ~$g$ is essential
and it is crucial to impose that the size of the g-balls is
preserved by Hamiltonian flow which leads to the only choice~$ K=
C|2^q I_q| h$ thanks to the properties of the solution of the
associate Hamilton Jacobi equation.

The key point in the generalization of the bilinear estimate is
the proof that for solutions of a variable coefficients wave
equation, microlocalization properties propagate nicely along the
Hamiltonian flows related to the wave operator; this point
follows from the choice of the metric used to localize in the
cotangent space of~$\RR^d$.

Finally to end the proof of the microlocal bilinear estimate, the
strategy consists to decompose the Cauchy data using unity
partition whose elements are supported in g-balls and then to
apply the product and the propagation theorems to concentrate on
real interaction (see the proof in the constant coefficient
case). Because of the fact that interaction in the product and
propagation of microlocalization are badly related, we need at
this step to recourse to a second microlocalization, which means
that we have to decompose again the interval on which we work on
sub intervals where the Hamiltonian flow is nearly constant .

\label{lastpage}

\end{document}